\documentclass[a4paper,11pt]{article}
\usepackage{stmaryrd}
\usepackage{amsfonts}
\usepackage{bbm}
\usepackage{float}
\usepackage{amscd}
\usepackage{mathrsfs}
\usepackage{latexsym,amssymb,amsmath,amscd,amscd,amsthm,amsxtra,xypic}
\usepackage[dvips]{graphicx}
\usepackage[utf8]{inputenc}
\usepackage[T1]{fontenc}
\usepackage{lmodern}
\usepackage{amssymb}
\usepackage[all]{xy}
\usepackage{nicefrac,mathtools,enumitem}
\usepackage{microtype}
\usepackage{lipsum}
\usepackage{mathtools}
\usepackage{xfrac}

\usepackage{geometry}
\usepackage{hyperref}

\usepackage{times}
\usepackage{amsthm}
\usepackage{amsmath}
\usepackage{amsfonts}
\usepackage{amssymb}
\usepackage{mathdots}
\usepackage{color}
\usepackage{mathrsfs}
\usepackage{stmaryrd}
\SetSymbolFont{stmry}{bold}{U}{stmry}{m}{n}

\usepackage{bm}

\usepackage{tikz-cd}
\usepackage{tikz}
\usetikzlibrary{matrix,arrows,decorations.pathmorphing}
\usepackage[all]{xy}
\usepackage{graphicx}
\usepackage{subfigure}
\usepackage[toc,page]{appendix}
\usepackage[usestackEOL]{stackengine}

\usepackage{theoremref}
\newtheorem{thm}{Theorem}[section]
\newtheorem{lem}[thm]{Lemma}
\newtheorem{cor}[thm]{Corollary}
\newtheorem{pro}[thm]{Proposition}

\theoremstyle{definition}

\newtheorem{rmk}[thm]{Remark}
\newtheorem{defi}[thm]{Definition}

\newcommand{\I}{{\mathrm{i}}}

\newcommand{\be }{\begin{equation}}
\newcommand{\ee }{\end{equation}}

\newcommand{\pf}{\noindent{\bf Proof.}\ }

\newcommand{\h}{\mathfrak h}

\newcommand{\U}{{\rm U}}

\def\qed{\hfill ~\vrule height6pt width6pt depth0pt}

\newcommand{\br}[1]{   [ \cdot,    \cdot  ]   }

\newcommand{\Ad}{\mathrm{Ad}}

\newcommand{\Herm}{{\rm Herm}}
\newcommand{\Sym}{{\rm Sym}}
\title{The Alekseev-Meinrenken diffeomorphism arising from the Stokes phenomenon}
\author{Xiaomeng Xu}
\date{}

\newcommand{\Addresses}{{
  \bigskip
  \footnotesize
\noindent \textsc{
School of Mathematical Sciences \& Beijing International Center
for Mathematical Research, Peking University, Beijing 100871, China}\par\nopagebreak
  \textit{E-mail address}: \texttt{xxu@bicmr.pku.edu.cn}
}}
\begin{document}

\maketitle
    
\begin{abstract}
The Alekseev-Meinrenken diffeomorphism is a distinguished diffeomorphism from the space of $n\times n$ Hermitian matrices to the space of $n\times n$ positive definite Hermitian matrices.
This paper derives the explicit expression of the diffeomorphism, via the Stokes phenomenon of meromorphic linear systems of ordinary differential equations with Poncar\'{e} rank $1$. 
\end{abstract}

\section{Introduction and the main result}
Let us consider the Lie algebra ${ \frak u}(n)$ of the unitary group ${\rm U}(n)$. The Ginzburg-Weinstein linearization theorem \cite{GW} states that the dual Poisson Lie group ${\rm U}(n)^*$, of the standard Poisson Lie group ${\rm U}(n)$ (see e.g. \cite{LW}), is Poisson isomorphic to the dual of the Lie algebra $\frak u(n)^*$. 
There are many proofs of Ginzburg-Weinstein linearization theorem, from the different perspectives of the cohomology calculation, Moser's trick in symplectic geometry, Stokes phenomenon, the quantum algebras and so on, see e.g., \cite{GW, Anton,AM, Boalch,BoalchG,EEM}. In particular, Alekseev and Meinrenken constructed the linearization via the Gelfand-Tsetlin integrable systems, Boalch's construction relies on the theory of meromorphic linear ODEs, while Enriquez-Etingof-Marshall Construction used the theory of quasi-Hopf algebras. 

It is natural to ask if these seemingly rather different methods are related to each other in some ways. We have been working on this problem with the belief that the connection between the different methods will bring new insights into both subjects. 

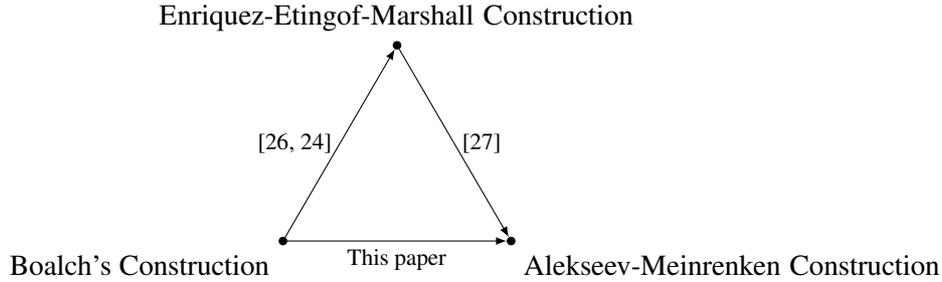
\[\begin{tikzpicture}[>=latex,mydot/.style={draw,circle,inner
    sep=1pt},every label/.style={scale=1},scale=1]

  \foreach \i in {0}{
  \node[mydot, fill=black, label=240:{Boalch's Construction}]                at (-1.5+9*\i, -2.595)    (p\i0) {};
  \node[mydot,fill=black,label=90:{Enriquez-Etingof-Marshall Construction}]      at (+9*\i,0) (p\i1) {};
  \node[mydot,fill=black,label=-60:{Alekseev-Meinrenken Construction}]     at (1.5+9*\i,-2.595)    (p\i2) {};
}
\begin{scope}[<-]
    \draw (p01)--node[left,scale=.8]{ \cite{Xu0, TLXu}} (p00);
    \draw (p02)--node[below,scale=.8]{This paper } (p00);
    \draw (p02)--node[right,scale=.8]{\cite{Xu00}} (p01);
\end{scope}
\end{tikzpicture}\]
Indeed, the pursuit of their connections (see \cite{Xu0, Xu00, TLXu}) leads to new realization of Gelfand-Tsetlin basis, and crystal basis in the representation of quantum groups via the Stokes phenomenon \cite{Xu, Xu1}, leads to the introduction of quantum Stokes matrices at arbitrary order poles and the quantization of the irregular Riemann-Hilbert maps \cite{Xuquantum}, and more importantly leads to answers to various analysis problems in the study of nonlinear isomonodromy deformation equations and WKB analysis \cite{Xu1, Xu2, ANXZ}. 

In this paper, we make one step further to build a precise connection between the Alekseev-Meinrenken construction and the Boalch construction. 

First, let us recall the Alekseev-Meinrenken construction. 
Let us identify the Lie algebra ${ \frak u}(n)$, consisting of skew-Hermitian matrices, with the space $\Herm(n)$ of $n\times n$ Hermitian matrices via the pairing $\langle A,\xi\rangle=2{\rm Im}({\rm tr}A\xi)$. Thus $\Herm(n)$ inherits a Poisson structure from the canonical linear (Kostant-Kirillov-Souriau) Poisson structure on ${\frak u}(n)^*$. Furthermore, the dual Poisson Lie group ${\rm U}(n)^*$, which is the group of complex
upper triangular matrices with strictly positive diagonal entries,
is identified with the space $\Herm^+(n)$ of positive definite Hermitian $n\times n$-matrices, by
taking the upper triangular matrix $X\in {\rm U}(n)^*$ to the positive Hermitian matrix $(X^*X)^{1/2}\in
\Herm^+(n)$. The Poisson manifolds ${\frak u}(n)^*\cong\Herm(n)$ and ${\rm U}(n)^*$ carry extra structures: Guillemin-Sternberg \cite{GS} introduced the Gelfand–Tsetlin integrable system on $\frak u(n)^*$; Flaschka-Ratiu \cite{FR} described a multiplicative Gelfand-Tsetlin system for the dual Poisson Lie group ${\rm U}(n)^*$. 
Later on in \cite{AM}, Alekseev and Meinrenken used these integrable systems to construct a distinguished Ginzburg-Weinstein linearization $\Gamma_{AM}$.

\begin{thm}\cite{AM}\label{11}
There exists a unique diffeomorphism
\[\Gamma_{AM}\colon\Herm(n)\cong\frak u(n)^*\rightarrow \Herm^+(n)\cong{\rm U}(n)^*,\] which intertwines the Gelfand-Tsetlin systems on both sides (and has one extra property, see Section \ref{AMdiff}). In particular, the map $\Gamma_{AM}$ is a Poisson isomorphism.
\end{thm}

Second, let us recall the Boalch construction.
Let $\h(\mathbb{R})$ (resp. $\h_{\rm reg}(\mathbb{R})$) denote the set of diagonal matrices with (resp. distinct) real eigenvalues. Let us consider the meromorphic linear system for a function $F(z)\in {\rm GL}_n(\mathbb{C})$,
\begin{eqnarray}\label{Stokeseq}
\frac{dF}{dz}=\left(\I  u-\frac{1}{2\pi\I }\frac{A}{z}\right)\cdot F,
\end{eqnarray}
where $\I =\sqrt{-1}$, $u\in\h(\mathbb{R})$ and $A\in{\Herm}(n)$. The system has an order two pole at $\infty$ and (if $A\neq 0$) a first order pole at $0$. It has a unique formal solution $\hat{F}(z)$ around
$z = \infty$. Then the standard theory of resummation states that there exist certain sectorial regions around $z=\infty$, such that on each of these sectors there is a unique (therefore canonical) holomorphic solution with the prescribed asymptotics $\hat{F}(z)$. These solutions are in general different (that reflects the Stokes phenomenon), and the transition between them can be measured by a pair of Stokes matrices $S_\pm(u,A)\in{\rm GL}_n(\mathbb{C})$. The Stokes matrices $S_+$ and $S_-$ are upper and lower triangular ones, and, due to the real condition $u\in \h(\mathbb{R})$ and $A\in\Herm(n)$, are complex conjugate to each other. See Section \ref{beginsection} for more details. Then the remarkable theorem of Boalch states that 

\begin{thm}\label{Boalchthm}\cite{Boalch}
For any fixed $u\in\h_{\rm reg}(\mathbb{R})$, the Riemann-Hilbert map
\begin{eqnarray}\label{introRH}
\nu(u)\colon\Herm(n)\cong\frak u(n)^*\rightarrow \Herm^+(n)\cong{\rm U}(n)^*; \ A\mapsto \nu(u,A):=S_-(u,A) S_+(u,A),
\end{eqnarray}
is a Poisson isomorphism (here since the Stokes matrices satisfy $S_-(u,A)=S_+(u,A)^\dagger$, the product $S_-S_+$ is a positive definite Hermitian matrix).
\end{thm}

In order to understand the relation between the Poisson diffeomorphisms $\Gamma_{AM}$ and $\nu(u)$, first note that $\Gamma_{AM}$ is a diffeomorphism compatible with the Gelfand-Tsetlin chain of subalebras 
\begin{equation}\label{GZchain}
\frak u(1)\subset \frak u(2)\subset \cdots \subset \frak u(n),
\end{equation}
where $\frak u(k)\subset \frak u(k+1)$ is the upper left corner, see Section \ref{appPoisson} for more details. While $\nu(u)$ is a family of diffeomorphisms depending on extra parameter space $u\in \h_{\rm reg}(\mathbb{R})$. Therefore, to relate $\nu(u)$ to $\Gamma_{AM}$, we need to break the symmetry of $\nu(u)$, by choosing specific $u$ such that the particular chain structure \eqref{GZchain} can come out. The naive observation is that if $u_n\gg u_{n-1}\gg \cdots \gg u_1$ in $u={\rm diag}(u_1,...,u_n)$, the $n \times n$ linear system \eqref{Stokeseq} seems to be decoupled into $n$ systems of rank $n,n-1, n-2, ..., 1$, according to the chain \eqref{GZchain}. In \cite{Xu1}, we made this naive observation a precise statement about the regularized limit of the Stokes matrices $S_\pm(u,A)$ as $u_n\gg u_{n-1}\gg \cdots \gg u_1$. The statement relies on the study of the nonlinear isomonodromy deformation equations of the linear system \eqref{Stokeseq}.

In this paper, motivated by the analysis results in \cite{Xu1}, we construct a family of explicit Poisson diffeomorphisms $\psi(u):\Herm(n)\rightarrow \Herm(n)$ parameterized by $u\in \h_{\rm reg}(\mathbb{R})$, via the Thimm action in the Gelfand-Tsetlin system. See Section \ref{endsection} for the explicit expression of $\psi(u)$. Then we show that the diffeomorphisms are such that the composed map $\Gamma_{AM}\circ \psi(u): \Herm(n)\rightarrow \Herm^+(n)$ is the first order approximation of the Riemann-Hilbert map $\nu(u)$ as $u_n\gg u_{n-1}\gg \cdots \gg u_1$. That is 
\begin{thm}\label{mainthm}
For any fixed $A\in\Herm(n)$, the Hermitian matrix function
\begin{equation}\label{AMnu}
\Gamma_{AM}\left(\psi(u,A)\right)=\nu(u,A)
+\sum_{k=2}^{n-1}\mathcal{O}\left(\frac{u_{k}-u_{k-1}}{u_{k+1}-u_{k}}\right),
\end{equation}
as $\frac{u_{k+1}-u_{k}}{u_{k}-u_{k-1}}\rightarrow +\infty$ for all $k=2,...,n-1$. Here we use the big $\mathcal{O}$ notation: each $\mathcal{O}\left(\frac{u_{k}-u_{k-1}}{u_{k+1}-u_{k}}\right)$ stands for a remainder whose norm is less than $M \times\left(\frac{u_{k+1}-u_{k}}{u_{k}-u_{k-1}}\right)^{-1}$ for a positive real number $M$ as $\frac{u_{k+1}-u_{k}}{u_{k}-u_{k-1}}$ big enough.
\end{thm}

As an application of Theorem \ref{mainthm}, we find the explicit formula of $\Gamma_{AM}$. See Section \ref{AMdiff} for more details. We denote by $\lambda^{(k)}_1\ge \lambda^{(k)}_2\ge\cdots\ge \lambda^{(k)}_k$ the ordered eigenvalues of the left-top $k\times k$ submatrix of $A$, and denote by $A_n={\rm diag}(\lambda^{(n)}_1,...,\lambda^{(n)}_n)$.

\begin{thm}\label{exAM}
The Alekseev-Meinrenken diffeomorphism is given by
\[\Gamma_{AM}\colon {\rm Herm}(n)\rightarrow {\rm Herm}^+(n)~; ~ A\mapsto\psi(A)e^{A_{n}}\psi(A)^{-1} \]
where $\psi$ is the pointwise multiplication $\psi(A)=\psi^{(1)}(A)\cdot\cdot\cdot \psi^{(n-1)}(A)\psi^{(n)}(A)$ of unitary matrices $\psi^{(k)}(A)$, and each map $\psi^{(k)}:\Herm(n)\rightarrow {\rm U}(n)$ is given by
\begin{align*}
\psi^{(k)}(A)_{ij} =&e^{\frac{\lambda^{(k-1)}_i-\lambda^{(k)}_j}{4}}\sqrt{\frac{ \prod_{v=1,v\ne i}^{k-1}{\rm sinh}\Big(\frac{\lambda^{(k-1)}_v-\lambda^{(k)}_j}{2}\Big)\prod_{v=1,v\ne j}^{k}{\rm sinh}\Big(\frac{\lambda^{(k-1)}_i-\lambda^{(k)}_v}{2}\Big)}{\prod_{v=1, v\ne j}^{k}{\rm sinh}\Big(\frac{\lambda^{(k)}_v-\lambda^{(k)}_j}{2}\Big)\prod_{v=1, v\ne i}^{k-1}{\rm sinh}\Big(\frac{\lambda^{(k-1)}_i-\lambda^{(k-1)}_v}{2}\Big)}}\\
&\times\frac{(-1)^{k-1+i}\Delta^{1,...,k-1}_{1,...,k-2,k}\left(A-\lambda^{(k-1)}_i\right)}{\sqrt{-{\prod_{l=1}^{k}(\lambda^{(k-1)}_i-\lambda^{(k)}_l)\prod_{l=1}^{k-2}(\lambda^{(k-1)}_i-\lambda^{(k-2)}_l)}}}, 
\end{align*}
for $1\le i\le k-1, 1\le j\le k$, and
\begin{align*} \psi^{(k)}(A)_{kj}&=e^{\frac{\lambda^{(k)}_j-\lambda^{(k-1)}_k}{4}}\sqrt{\frac{\prod_{v=1}^{k-1}{\rm sinh}\Big(\frac{\lambda^{(k-1)}_v-\lambda^{(k)}_j}{2}\Big)}{\prod_{v=1, v\ne j}^{k}{\rm sinh}\Big(\frac{\lambda^{(k)}_v-\lambda^{(k)}_j}{2}\Big)}}, \ \ \ \text{ for } \ 1\le j\le k,\\
\psi^{(k)}(A)_{ii}&=1, \ \ \ \ \ \ \text{for} \ k<i\le n, \ \ \ \ \text{ and }
\psi^{(k)}(A)_{ij}=0, \ \ \ \ \ \ \text{otherwise}.
\end{align*}
Here $\Delta^{1,...,k-2,k-1}_{1,...,k-2,k}(A-{\lambda^{(k-1)}_i})$ is the $(k-1)\times (k-1)$ minor of the matrix $A-{\lambda^{(k-1)}_i}\cdot {\rm Id}_{n}$ formed by the first $k-1$ rows and $1,...,k-2,k$ columns (${\rm Id}_n$ is the rank $n$ identity matrix).
\end{thm}

The left hand side of the identity \eqref{AMnu} is characterized by the algebraic integrable system, while the right hand side is highly transcendental. Therefore, Theorem \ref{mainthm} provides us the possibility to study the Stokes phenomenon, as well as the associated isomonodromy deformation equation, via algebraic methods. For example, in the above discussion, we have chosen a particular Gelfand-Tsetlin chain of subalgebras. Instead, we can take another chain, the chain of subalgebras $\frak u(1)\subset \frak u(2)\subset \cdots \subset \frak u(n),$
where $\frak u(k)\subset \frak u(k+1)$ is the lower right (instead of upper left) corner. Accordingly, we can define Alekseev-Meinrenken diffeomorphism from $\Herm(n)$ to $\Herm^+(n)$ compatible with the corresponding Gelfand-Tsetlin systems. Then such diffeomorphism is related to the Riemann-Hilbert map $\nu(u)$ as $u_1\gg u_2\gg \cdots \gg u_1$ in the same way as Theorem \ref{mainthm}. The comparison of the results, associated to the two different choices of Gelfand-Tsetlin chains, can be used to characterize the connection formula of the nonlinear isomonodromy deformation equations between the two different asymptotic zones $\frac{u_{k+1}-u_{k}}{u_{k}-u_{k-1}}\rightarrow \infty$ for $k=2,...,n-1$ and $\frac{u_{k+1}-u_{k}}{u_{k}-u_{k-1}}\rightarrow 0$ for $k=2,...,n-1$. See \cite{Xu} for more details.

Another application of Theorem \ref{exAM} is that it enables us to write down an explicit Moser flow of the Hamiltonian vector fields on the coadjoint orbit of $\frak{u}(n)^*$ constructed by Alekseev in \cite{Anton}. Therefore, it will relate Alekseev's construction using Moser's trick \cite{Anton} to the other constructions of Ginzburg-Weinstein linearization, and thus bring more insights into different subjects.
We also remark that it is interesting to study the behaviour of the expression in Theorem \ref{exAM} under the toric degeneration of Gelfand-Tsetlin systems, see e.g., \cite{NNU}. We leave the study of the relations between the Stokes phenomenon, Moser flow and toric degeneration to a future work.

\vspace{2mm}
The organization of the paper is as follows. Section \ref{beginsection} gives the preliminaries of Stokes data of meromorphic linear systems, and recalls the expression of Stokes matrices via the boundary value of the associated nonlinear isomonodromy equation. Section \ref{appPoisson} recalls the Gelfand-Tsetlin systems, and gives a proof of Theorem \ref{mainthm}, as well as Theorem \ref{exAM}. 

\subsection*{Acknowledgements}
\noindent
We would like to thank Anton Alekseev and Eckhard Meinrenken for their discussion and encouragement. The author is supported by the National Key Research and Development Program of China (No. 2021YFA1002000) and by the National Natural Science Foundation of China (No. 12171006).

\section{Stokes phenomenon and monodromy data}\label{beginsection}
In Section \ref{Canonicalsol}, we recall the canonical solutions of the equation \eqref{Stokeseq}. In Section \ref{defiStokes}, we introduce the Stokes matrices and connection matrices of the linear systems, as well as the monodromy relation relating connection matrices to Stokes matrices. In Sections \ref{sec:iso} and \ref{firstsec}, we recall the boundary value of the solutions of the isomonodromy deformation equation of \eqref{Stokeseq}, and express the Stokes matrices of \eqref{Stokeseq} via the boundary value of the associated isomonodromy equation.
\subsection{Canonical solutions}\label{Canonicalsol}
Let $\h(\mathbb{R})$ (resp. $\h_{\rm reg}(\mathbb{R})$) denote the set of diagonal matrices with (resp. distinct) real eigenvalues. Let us consider the meromorphic linear system \eqref{Stokeseq}.

\begin{defi}\label{Stokesrays}
The {\it Stokes supersectors} of the system
are the two sectors ${\rm Sect}_+:=\{z\in\mathbb{C}~|~-\pi<{\rm arg}(z)<\pi\}$ and ${\rm Sect}_-=\{z\in\mathbb{C}~|~-2\pi<{\rm arg}(z)<0\}$.
\end{defi}

Let us choose the branch of ${\rm log}(z)$, which is real on the positive real axis, with a cut along the nonnegative imaginary axis $\I  \mathbb{R}_{\ge 0}$. Then by convention, ${\rm log}(z)$ has imaginary part $-\pi$ on the negative real axis in ${\rm Sect}_-$.
\begin{thm}\label{uniformresum}
For any $u\in\h_{\rm reg}(\mathbb{R})$, on ${\rm Sect}_\pm$ there
is a unique fundamental solution $F_\pm:{\rm Sect}_\pm\to {\rm GL}(n,\mathbb{C})$ of equation \eqref{Stokeseq} such that 
\begin{align*}
\lim_{z\rightarrow\infty}F_+(z;u)\cdot e^{-\I  uz}\cdot z^{\frac{[A]}{2\pi\I }}&={\rm Id}_n, \ \ \ \text{as} \ \ \ z\in {\rm Sect}_+,
\\
\lim_{z\rightarrow\infty}F_-(z;u)\cdot e^{-\I  uz}\cdot z^{\frac{[A]}{2\pi\I }}&={\rm Id}_n, \ \ \ \text{as} \ \ \ z\in {\rm Sect}_-,
\end{align*}
Here ${\rm Id}_n$ is the rank n identity matrix, and $[A]$ is the diagonal part of $A$.
\end{thm}

\subsection{Stokes matrices}\label{secSC}
For any $\sigma\in S_n$, let us denote by $U_{\sigma}$ the component $\{(u_1, . . . , u_n)~|~ u_{\sigma(1)}<\cdot\cdot\cdot < u_{\sigma(n)}\}$ of $\h_{\rm reg}(\frak{sl}_n)(\mathbb{R})$, and denote by $P_\sigma\in{\rm GL}_n$ the corresponding permutation matrix.

\begin{defi}\label{defiStokes}
For any $u\in U_\sigma$, the {\it Stokes matrices} of the system \eqref{Stokeseq} (with respect
to ${\rm Sect}_+$ and the chosen branch of ${\rm log}(z)$) are the elements $S_\pm(u,A)\in {\rm GL}(n)$ determined by
\[F_+(z)=F_-(z)\cdot e^{-\frac{[A]}{2}}P_\sigma S_+(u,A)P_\sigma^{-1}, \  \ \ \ \ 
F_-(ze^{-2\pi \I })=F_+(z)\cdot P_\sigma S_-(u,A)P_\sigma^{-1}e^{\frac{[A]}{2}},
\]
where the first (resp. second) identity is understood to hold in ${\rm Sect}_-$
(resp. ${\rm Sect}_+$) after $ F_+$ (resp. $F_-$)
has been analytically continued anticlockwise around $z=\infty$. 
\end{defi} 
The prescribed asymptotics of $F_\pm(z)$ at $z=\infty$, as well as the identities in Definition \ref{defiStokes}, ensures that the Stokes matrices $S_+(u,A)$ and $S_-(u,A)$ are upper and lower triangular matrices respectively. see e.g., \cite[Chapter 9.1]{Wasow}. Furthermore, the following lemma follows from the fact that if $F(z)$ is a solution, so is $F(\bar{z})^\dagger$, see \cite{Boalch}.
\begin{lem}
Let $S_+(u,A)^\dagger$ denote the conjugation transpose of $S_+(u,A)$, then $S_-(u,A)=S_+(u,A)^\dagger$.
\end{lem}

Since the system \eqref{Stokeseq} is non-resonant, i.e., no two eigenvalues of $\frac{A}{2\pi\I }$ for $A\in\Herm(n)$ are differed by a positive integer, we have (see e.g \cite[Chapter 2]{Wasow}).
\begin{lem}\label{le:nr dkz}
There is a unique holomorphic fundamental solution
$F_0(z;u,A)\in {\rm GL}(n)$ of the system \eqref{Stokeseq} on a neighbourhood of $\infty$ slit along $\I  \mathbb{R}_{\ge 0}$, such that $F_0\cdot z^{\frac{A}{2\pi\I }}\rightarrow {\rm Id}_n$ as $z\rightarrow 0$.
\end{lem}

\begin{defi}\label{connectionmatrix}
The {\it connection 
matrix} $C(u,A)\in {\rm GL}_n(\mathbb{C})$ of the system \eqref{Stokeseq} (with respect to ${\rm Sect}_+$) is determined by 
\[F_0(z;u,A)=F_+(z;u,A)\cdot C(u,A), \]
as $F_0(z;u,A)$ is
extended to the domain of definition of $F_+(z;u,A)$.
\end{defi}
The connection matrix $C(u,A)$ is valued in ${\rm U}(n)$ (see e.g., \cite[Lemma 29]{Boalch}).
Thus for any fixed $u$, by varying $A\in{\Herm}(n)$ we obtain the connection map
\begin{eqnarray}\label{Cmap}
C(\cdot, u)\colon{\rm Herm}(n)\rightarrow {\rm U}(n).\end{eqnarray}

In a global picture, the connection matrix is related to the Stokes matrices
by the following monodromy relation, which follows from the fact that a simple negative loop (i.e., in clockwise direction) around $0$ is a simple positive loop (i.e., in anticlockwise direction) around $\infty$: for any $u\in U_\sigma\subset \h_{\rm reg}(\mathbb{R})$,
\begin{eqnarray}\label{monodromyrelation}
C(u,A)e^{A}C(u,A)^{-1}=P_\sigma S_-(u,A)S_+(u,A) P_\sigma^{-1}.
\end{eqnarray}

\subsection{Isomonodromy deformation}\label{sec:iso}
In this subsection, we recall some facts about the theory of isomonodromy deformation. In general, the Stokes matrices $S_\pm(u,A)$ of the system \eqref{Stokeseq} will depend on the irregular term $u$.
The isomonodromy deformation (also known as monodromy preserving) problem is to find the matrix valued function $\Phi(u)$ such that the Stokes matrices $S_\pm(u,\Phi(u))$ are (locally) constant. In particular, the following definition and proposition are known. See more detailed discussions in e.g., \cite{JMMS, JMU, Dubrovin, BoalchG}.
\begin{defi}
The isomonodromy equation is the differential equation for a matrix valued function $\Phi(u): \h_{\rm reg}(\mathbb{R})\rightarrow {\rm Herm}(n)$
\begin{eqnarray}\label{isoeq}
\frac{\partial \Phi}{\partial u_k} =\frac{1}{2\pi\I }[\Phi,{\rm ad}^{-1}_u{\rm ad}_{E_{k}}\Phi], \ \text{for all} \ k=1,...,n.\end{eqnarray} 
Here $E_k$ is the $n\times n$ diagonal matrix whose $(k,k)$-entry is $1$ and other entries are $0$.
Note that ${\rm ad}_{E_{k}}\Phi$ takes values in the space ${\frak {gl}}_n^{od}$ of off diagonal matrices and that ${\rm ad}_u$
is invertible when restricted to ${\frak {gl}}_n^{od}$.
\end{defi}

Set $\Phi(u)=(\phi_{ij}(u))$, then in terms of the components, the equation \eqref{isoeq} becomes
\begin{eqnarray*}
  \frac{\partial}{\partial u_k} \phi_{i j} (u)
  & = & 
  \frac{1}{2\pi\I }\left(
  \frac{1}{u_k - u_j}-\frac{1}{u_k - u_i} \right) 
  \phi_{i k}(u) \phi_{k j}(u), \quad i, j \neq k,
  \label{uisoeqij}
  \\
  \frac{\partial}{\partial u_k} \phi_{i k} (u)
  & = & \frac{1}{2\pi\I }
  \sum_{j \neq k} 
  \frac{\phi_{i j}(u)\phi_{j k}(u)-\delta_{i j}\phi_{k k}(u) \phi_{j k}(u)}{u_k - u_j} ,
  \quad i \neq k ,
  \label{uisoeqik}
  \\
  \frac{\partial}{\partial u_k} \phi_{k j} (u)
  & = & \frac{1}{2\pi\I }
  \sum_{i \neq k} 
  \frac{\delta_{i j}
  \phi_{k k}(u)\phi_{k i}(u)-\phi_{k i}(u)\phi_{i j}(u)}{u_k - u_i},
  \quad j \neq k ,
  \label{uisoeqkj}\\
  \frac{\partial}{\partial u_k} \phi_{k k}(u) 
  & = & 
  0.
  \label{uisoeqkk}
\end{eqnarray*}
\begin{pro}\label{isomonodef}
For any solution $\Phi(u)$ of the isomonodromy equation, the Stokes matrices $S_\pm(u,\Phi(u))$
are locally constants (independent of $u$). 
\end{pro}
\begin{rmk}
Following Miwa \cite{Miwa}, the $\frak{gl}_n$-valued solutions $\Phi(u)$ of the equation \eqref{isoeq} with $u_1,...,u_n\in \mathbb{C}$ have the strong Painlev\'{e} property: they are multi-valued meromorphic functions of $u_1,...,u_n$ and the branching occurs when $u$ moves along a loop around the fat diagonal
\[\Delta=\{(u_1,...,u_n)\in \mathbb{C}^n~|~u_i = u_j, \text{for some } i\ne j \}.\]
Then, according to Boalch \cite{Boalch}, when restricts to the real case, the $\Herm(n)$-valued solutions $\Phi(u)$ of \eqref{isoeq} are real analytic on each connected component of $u\in \h_{\rm reg}(\mathbb{R})$.
\end{rmk}

\subsection{The boundary value and explicit Stokes matrices}\label{firstsec}
Now let us consider the $n\times n$ system of partial differential equations for a function $F(z,u)\in GL(n)$ 
\begin{eqnarray}\label{introisoStokeseq1}
\frac{\partial F}{\partial z}&=&\left(\I  u-\frac{1}{2\pi\I }\frac{\Phi(u)}{z}\right)\cdot F,\\
\label{introisoStokeseq2}
\frac{\partial F}{\partial u_k}&=&\left(\I  E_kz-\frac{1}{2\pi\I } {\rm ad}^{-1}_u{\rm ad}_{E_{k}}\Phi(u)\right)\cdot F, \ \text{for all} \ k=1,...,n.
\end{eqnarray}
where the residue $\Phi(u)\in {\rm Herm}(n)$ is a solution of the isomonodromy equation \eqref{isoeq}. One checks that \eqref{isoeq} is the compatibility condition of the above PDE system. The boundary value and monodromy problem of the system were studied in \cite{Xu}.
\begin{thm}\cite[Theorem 1.1]{Xu}\label{isomonopro}
For any $\Herm(n)$-valued solution $\Phi(u)$ of the isomonodromy equation \eqref{isoeq} on the connected component $U_{\rm id}:=\{u\in \h_{\rm reg}(\mathbb{R})~|~u_1<\cdots <u_n\}$, there exists a unique constant $\Phi_0\in\Herm(n)$ such that as the real numbers $\frac{u_{k+1}-u_{k}}{u_{k}-u_{k-1}}\rightarrow +\infty$ for all $k=2,...,n-1$,
\begin{align}\label{firstasy}
\Phi(u)=&{\rm Ad}{\left((u_2-u_1)^{\frac{\delta_1(\Phi_0)}{2\pi\I }}\cdot \overrightarrow{\underset{k=2,...,n-1}{\prod} }\left(\frac{u_{k+1}-u_{k}}{u_{k}-u_{k-1}}\right)^{\frac{\delta_k(\Phi_0)}{2\pi\I }}\right)}\Phi_0 +O\left(\left(\frac{u_{k+1}-u_{k}}{u_{k}-u_{k-1}}\right)^{-1};k=2,...,n-1\right),
\end{align}
where ${\rm Ad}(g)X=gXg^{-1}$ for any $g\in U(n)$ and $X\in \Herm(n)$, the product $\overrightarrow{\prod}$ is taken with the index $i$ to the right of $j$ if $i>j$. And $\delta_k(\Phi)$ is the Hermitian matrix with entries
\[\delta_k(\Phi)_{ij}=\left\{
          \begin{array}{lr}
             \Phi_{ij},   & \text{if} \ \ 1\le i, j\le k, \ \text{or} \ i=j  \\
           0, & \text{otherwise}.
             \end{array}
\right.\] Furthermore, given any $\Phi_0\in\Herm(n)$ there exists a unique real analytic solution $\Phi(u)$ of \eqref{isoeq} with the prescribed asymptotics \eqref{firstasy}.
\end{thm}
Therefore, $\Phi_0\in\Herm(n)$ parameterizes the Hermitian matrix valued solutions of \eqref{isoeq} on $U_{\rm id}$. 
We then denote by $\Phi(u;\Phi_0)$ the solution of \eqref{isoeq} with the prescribed asymptotics $\Phi_0$ in the sense of Theorem \ref{isomonopro}. 

In \cite{Xu}, the Stokes matrices of the systems \eqref{introisoStokeseq1}-\eqref{introisoStokeseq2} are given explicitly in terms of $\Phi_0$. Let us denote by $\{\lambda^{(k)}_i\}_{i=1,...,k}$ the eigenvalues of the left-top $k\times k$ submatrix of $\Phi_0$, and $(\Phi_0)_{k+1,k+1}$ the $k+1$-th diagonal element. First, it was proved that
\begin{thm}\cite[Theorem 3.14]{Xu}\label{relautomorphism}
For any solution $\Phi(u;\Phi_0)$ of the isomonodromy equation \eqref{isoeq} on $u\in U_{\rm id}$ with the prescribed asymptotics $\Phi_0\in\Herm(n)$ (as in Theorem \ref{isomonopro}), we have
\[S_{-}(u, \Phi(u;\Phi_0)) S_{+}(u, \Phi(u;\Phi_0))=\left(\overrightarrow{\prod_{k=1,...,n}}C\left(E_k,\delta_k(\Phi_0)\right)\right)\cdot e^{\Phi_0}\cdot \left(\overrightarrow{\prod_{k=1,...,n}}C\left(E_k,\delta_k(\Phi_0)\right)\right)^{-1},\] 
where the product $\overrightarrow{\prod}$ is taken with the index $i$ to the right of $j$ if $j<i$, and for any $k=1,...,n$, $C\left(E_k,\delta_k(\Phi_0)\right)\subset U(n)$ is the connection matrix of the $n\times n$ system
\begin{eqnarray}\label{eq:relativeStokes}
\frac{dF}{dz}=\left(\I  E_k-\frac{1}{2\pi \I }\frac{\delta_k(\Phi_0)}{z}\right)F,
\end{eqnarray}
where $E_k={\rm diag}(0,...,0,1,0,..,0)$ with $1$ at the $k$-th position. 
\end{thm} 

Since for each $k$ the equation \eqref{eq:relativeStokes} can be solved explicitly via the confluent hypergeometric function $_kF_k$, thus the connection matrices $C\left(E_k,\delta_k(\Phi_0)\right)$ can be computed explicitly via the known asymptotic formula of the confluent hypergeometric function $_kF_k$. Then together with a manipulation of Gauss decomposition, the formula in Theorem \ref{relautomorphism} leads to
\begin{thm}\cite[Theorem 1.3]{Xu}\label{Stokesexp}
The sub-diagonals of the Stokes matrices $S_\pm(u,\Phi(u;\Phi_0))$ are given by
\begin{align*}
(S_+)_{k,k+1}&=2\pi\I  e^{\frac{\small{(\Phi_0)_{kk}+(\Phi_0)_{k+1,k+1}}}{4}} \\
\times &\sum_{i=1}^k\frac{\prod_{l=1,l\ne i}^{k}\Gamma\left(1+\frac{\lambda^{(k)}_l-\lambda^{(k)}_i}{2\pi \I }\right)}{\prod_{l=1}^{k+1}\Gamma\left(1+\frac{\lambda^{(k+1)}_l-\lambda^{(k)}_i}{2\pi \I }\right)}\frac{\prod_{l=1,l\ne i}^{k}\Gamma\left(\frac{\lambda^{(k)}_l-\lambda^{(k)}_i}{2\pi \I }\right)}{\prod_{l=1}^{k-1}\Gamma\left(1+\frac{\lambda^{(k-1)}_l-\lambda^{(k)}_i}{2\pi \I }\right)}\cdot {(-1)^{k+i}\Delta^{1,...,k-1,k}_{1,...,k-1,k+1}\left(\frac{\Phi_0-\lambda^{(k)}_i}{2\pi\I }\right)},\\
(S_-)_{k+1,k}&=-2\pi\I  e^{\frac{\small{(\Phi_0)_{kk}+(\Phi_0)_{k+1,k+1}}}{4}}\\
\times &
\sum_{i=1}^k \frac{\prod_{l=1,l\ne i}^{k}\Gamma\left(1-\frac{\lambda^{(k)}_l-\lambda^{(k)}_i}{2\pi \I }\right)}{\prod_{l=1}^{k+1}\Gamma\left(1-\frac{\lambda^{(k+1)}_l-\lambda^{(k)}_i}{2\pi \I }\right)}\frac{\prod_{l=1,l\ne i}^{k}\Gamma\left(-\frac{\lambda^{(k)}_l-\lambda^{(k)}_i}{2\pi \I }\right)}{\prod_{l=1}^{k-1}\Gamma\left(1-\frac{\lambda^{(k-1)}_l-\lambda^{(k)}_i}{2\pi \I }\right)}\cdot {(-1)^{k+i}\Delta^{1,...,k-1,k+1}_{1,...,k-1,k}\left(\frac{\lambda^{(k)}_i-\Phi_0}{2\pi\I }\right)}.
\end{align*}
where $k=1,...,n-1$ and $\Delta^{1,...,k-1,k}_{1,...,k-1,k+1}\left(\frac{\Phi_0-{\lambda^{(k)}_i}}{2\pi\I  }\right)$ is the $k$ by $k$ minor of the matrix $\frac{\Phi_0-{\lambda^{(k)}_i}}{2\pi\I  }$ formed by the first $k$ rows and $1,...,k-1,k+1$ columns. 
Furthermore, the other entries are given by explicit expressions via the sub-diagonal ones.
\end{thm}

\begin{defi}\label{caterpillar}
For any $\Phi_0\in \Herm(n)$, we introduce $S_{\pm}(u_{\rm cat}, \Phi_0):=S_{\pm}(u, \Phi(u;A))$, called the Stokes matrices at a caterpillar point (with respect to the connected component $U_{\rm id}$).
\end{defi}
\begin{rmk}
Theorem \ref{isomonopro} and \ref{Stokesexp} solve the boundary and monodromy problems for the isomonodromy equation of the $n\times n$ meromorphic linear system of ordinary differential equations with Poncar\'{e} rank $1$. In the case $n=3$, they recover Jimbo's formula for Painlev\'e VI equation \cite{Jimbo}.
\end{rmk}
\begin{rmk}
More generally, the regularized limit of the Stokes matrices $S_\pm(u,A)$ of equation \eqref{Stokeseq}, as some components $u_i$ of $u={\rm diag}(u_1,...,u_n)$ collapse in a comparable speed, was studied in \cite{Xu}. 
The prescription of the regularized limits is controlled by the geometry of the De Concini-Procesi wonderful compactification space $\widetilde{\h_{\rm reg}}(\mathbb{R})$. And it completely describes the asymptotic behaviour of the Stokes matrices at the singularities $u_i=u_j$. Here given a finite set of subspaces of a vector space, the De Concini-Procesi space \cite{dCP} replaces the set of subspaces by a divisor with normal crossings, and leaves the complement of these subspaces unchanged. As for the hyperplanes $u_i=u_j$ for all indices $i\ne j$, 
the associated De Concini-Procesi space $\widetilde{\h_{\rm reg}}(\mathbb{R})$ contains $\h_{\rm reg}(\mathbb{R})$ as an open part, and 
roughly speaking, a point in the boundary $\widetilde{\h_{\rm reg}}(\mathbb{R})\setminus{\h_{\rm reg}}(\mathbb{R})$ is a limit point $u={\rm diag}(u_1,...,u_n)$, where some $u_i$ collapse in a comparable speed. In particular, the limit of $u={\rm diag}(u_1,...,u_n)$, as $\frac{u_{k+1}-u_{k}}{u_{k}-u_{k-1}}\rightarrow +\infty$ for all $k=2,...,n-1$ and $u_2-u_1\rightarrow 0$, is a point $u_{\rm cat}$ in the $0$-dimensional stratum of $\widetilde{\h_{\rm reg}}(\mathbb{R})$, called a caterpillar point (see \cite{Sp}, page 16). 

Note that $u$ can approach to a boundary point, for example $u_{\rm cat}$, from different connected components $U_\sigma$ of ${\h_{\rm reg}}(\mathbb{R})$. When study the regularized limits of Stokes matrices at a boundary point, for example $u_{\rm cat}$, we should specialize the connected component from which we take the limit. It is because that as $u$ approaches to a fixed boundary point from two different components, the regularized limits can be different, see \cite{Xu} for more details. It explains the reason why we stress that the Riemann-Hilbert map at $u_{\rm cat}$ in Definition \ref{RHBcat} is with respect to $U_{\rm id}$.
\end{rmk}

As a consequence, for large $\frac{u_{k+1}-u_{k}}{u_{k}-u_{k-1}}$, letting $A=\Phi(u;\Phi_0)$ in Theorem \ref{relautomorphism} and Theorem \ref{Stokesexp} leads to 
\begin{pro}\cite[Proposition 3.32]{Xu1}\label{leadingterm}
For any fixed $A\in\Herm(n)$, we have that as $u\in U_{\rm id}$ and $\frac{u_{k+1}-u_{k}}{u_{k}-u_{k-1}}\rightarrow+\infty$,
\begin{align}\label{relimitGZ}
S_\pm(u,A)=S_\pm\left(u_{\rm cat}, \hspace{2mm} g(u;A)\cdot A\cdot g(u;A)^{-1}\right)+\sum_{k=2}^{n-1}\mathcal{O}(\frac{u_{k}-u_{k-1}}{u_{k+1}-u_{k}}),
\end{align}
where the unitary matrix
\begin{equation}\label{gu}
g(u;A):=(\frac{1}{u_{2}-u_{1}})^{\frac{(\delta_1(A))}{2\pi\I }}\cdot\overrightarrow{\underset{k=2,...,n-1}{\prod} }(\frac{u_{k}-u_{k-1}}{u_{k+1}-u_{k}})^{\frac{(\delta_k(A))}{2\pi\I }}.
\end{equation}
\end{pro}
This proposition will be used in Section \ref{endsection} to find the explicit expression of the diffeomorphism $\psi(u)$ in Theorem \ref{mainthm}. 
\subsection{Explicit Riemann-Hilbert map at a caterpillar point}

Theorem \ref{isomonopro} gives a parameterization of the Hermitian matrix valued solutions of the isomonodromy equation \eqref{isoeq}, and Theorem \ref{mainthm} computes explicitly the Stokes matrices of the corresponding linear equation \ref{introisoStokeseq1} via the parameterization. Therefore, we obtain an explicit Riemann-Hilbert map (a diffeomorphism) from $\Herm(n)$ to the space of Stokes matrices via the equivalences
\begin{eqnarray} \nonumber
    \Big\{\Phi_0\in\Herm(n)\Big\}&\Longleftrightarrow &\Big\{\text{solutions $\Phi(u;\Phi_0)\in\Herm(n)$ of the isomonodromy equation \eqref{isoeq} on $U_{\rm id}$}\Big\}\\ \nonumber &\Longleftrightarrow&  \Big\{\text{linear systems of PDEs \eqref{introisoStokeseq1} and \eqref{introisoStokeseq2} with $u\in U_{\rm id}$}\Big\} \\ \nonumber \label{Poissonnature} & \Longleftrightarrow & \Big\{\text{space of Stokes matrices $S_\pm(u,\Phi(u;\Phi_0))$ with $u\in U_{\rm id}$}\Big\}. 
\end{eqnarray}
It motivates
\begin{defi}\cite{Xu}\label{RHBcat}
The Riemann-Hilbert map at a caterpillar point $u_{\rm cat}$ (with respect to $U_{\rm id}$) is defined by 
\[\nu(u_{\rm cat})\colon\Herm(n)\cong\frak u(n)^*\rightarrow \Herm^+(n)\cong{\rm U}(n)^*; \ A\mapsto S_{-}(u_{\rm cat},A) S_+(u_{\rm cat},A).\]

Then by Theorem \ref{relautomorphism}, we can also write the Riemann-Hilbert map via the connection matrices
\begin{equation}\label{nuC}
\nu(u_{\rm cat},A)=\left(\overrightarrow{\prod_{k=1,...,n}}C\left(E_k,\delta_k(A)\right)\right)\cdot e^{A}\cdot \left(\overrightarrow{\prod_{k=1,...,n}}C\left(E_k,\delta_k(A)\right)\right)^{-1}.
\end{equation}
\end{defi}
The identity \eqref{nuC} will be used in Section \ref{AMdiff} to find the explicit expression of $\Gamma_{AM}$.

\section{The explicit Alekseev-Meinrenken diffeomorphism}\label{appPoisson}
In this section, we prove Theorem \ref{exAM} and Theorem \ref{mainthm}. Section \ref{GTcoor} introduces the Gelfand-Tsetlin systems. Section \ref{rRHvsGZ} shows that the Riemann-Hilbert maps at caterpillar points intertwine Gelfand-Tsetlin systems and their multiplicative analogs. Then Section \ref{AMdiff} uses the phase transformation of Gelfand-Tsetlin systems to derive the explicit formula of the map $\Gamma_{AM}$. In the end, Section \ref{endsection} proves Theorem \ref{mainthm}.

\subsection{The Gelfand-Tsetlin coordinates}\label{GTcoor}
{\bf Gelfand-Tsetlin maps.} For $k\le n$ let $A^{(k)}\in \Herm(k)$ denote the upper left $k$-th submatrix (upper left $k\times k$ corner) of a Hermitian matrix $A\in \Herm(n)$,
and $\lambda^{(k)}_i(A)$-its ordered set of
eigenvalues, $\lambda_1^{(k)}(A)\ge \cdots\ge \lambda_k^{(k)}(A)$. 
The map 
\begin{equation}\label{eq:momentmap}
\lambda\colon \Herm(n)\to \mathbb{R}^{\frac{n(n+1)}{2}},
\end{equation}
taking $A$ to the collection of numbers $\lambda_i^{(k)}(A)$ for $1\le i\le
k\le n$, is continuous and is called the Gelfand-Tsetlin map.
Its image $\mathcal{C}(n)$ is the Gelfand-Tsetlin cone, cut out by the following inequalities,
\begin{equation}\label{eq:cone}
\lambda_i^{(k+1)}\ge \lambda_i^{(k)}\ge \lambda_{i+1}^{(k+1)},\ \ 1\le
i\le k\le n-1.
\end{equation}

\vspace{2mm}
{\bf Thimm torus actions.}
Let $\mathcal{C}_0(n)\subset \mathcal{C}(n)$ denote the subset where all of the eigenvalue inequalities
\eqref{eq:cone} are strict. Let $\Herm_0(n):=\lambda^{-1}(\mathcal{C}_0(n))$ be the corresponding dense open subset of $\Herm(n)$. The $k$-torus 
$T(k)\subset {\rm U}(k)$ of diagonal matrices acts on 
$\Herm_0(n)$ as follows,  
\begin{equation}\label{eq:taction}
 t\bullet A=\Ad_{U^{-1} t U}A,\ \ \ \ t\in T(k),\ A\in \Herm_0(n).
\end{equation}
Here $U\in {\rm U}(k)\subset {\rm U}(n)$ is a unitary matrix such that $\Ad_{U}A^{(k)}$ is
diagonal, with entries $\lambda_1^{(k)},\ldots,\lambda_k^{(k)}$. The action
is well-defined since $U^{-1}t U$ does not depend on the choice of
$U$, and preserves the Gelfand-Tsetlin map \eqref{eq:momentmap}.  The actions
of the various $T(k)$'s commute, hence they define an action of the Gelfand-Tsetlin torus $T(1)\times \cdots \times T(n-1)\cong \U(1)^\frac{(n-1)n}{2}.$
Here the torus $T(n)$ is excluded, since the 
action \eqref{eq:taction} is trivial for $k=n$.  

\vspace{2mm}
{\bf Action-angle coordinates.}
If $A\in \rm Herm_0(n)$, then there exists a unique unitary matrix $P_k(A)\in {\rm U}(k)\subset {\rm U}(n)$, whose entries in the $k$-th row are positive and real, such that the upper left $k$-th submatrix of $A_k:=P_k(A)^{-1}AP_k(A)$ is the diagonal matrix ${\rm diag}(\lambda^{(k)}_1,...,\lambda^{(k)}_k)$, i.e., 
\begin{eqnarray}\label{zhuazi}
    A_k=P_k(A)^{-1}AP_k(A)=\begin{pmatrix}
    \lambda^{(k)}_1 & & & &a^{(k)}_1 & \cdots\\
    & \ddots & && \vdots & \cdots\\
    & &  \lambda^{(k)}_k && a^{(k)}_k & \cdots \\
    &&&&\\
    \overline{a_1^{(k)}}&\cdots &\overline{a_k^{(k)}}&& \lambda_{k+1}^{(k)} & \cdots\\
    \cdots & \cdots & \cdots &&\cdots & \cdots
    \end{pmatrix}.
\end{eqnarray}
The $(i,k+1)$ entries $a^{(k)}_i(A)$, for $1\le i\le k\le n-1$, are seen as functions on $\Herm_0(n)$.
\begin{defi}\label{GZfunctions}
The functions $\{\lambda^{(k)}_i\}_{1\le i\le k\le n}$ and $\{\psi^{(k)}_i={\rm Arg}(a^{(k)}_i)\}_{1\le i\le k\le n-1}$ on ${\rm Herm}_0(n)$ are called the Gelfand-Tsetlin action and angle coordinates.
\end{defi}

For any $2\le k\le n$, the $n\times n$ matrix $P_k(A)$ has the form
\begin{eqnarray}\label{clsP}
&&(P_k)_{ij}:=\frac{{(-1)^{i+j}}\Delta_{1,...,\hat{i},...,k}^{1,...,k-1}\left(A-\lambda^{(k)}_j\right)}{\sqrt{\prod_{l=1,l\ne i}^k(\lambda^{(k)}_j-\lambda^{(k)}_l)\prod_{l=1}^{k-1}(\lambda^{(k)}_j-\lambda^{(k-1)}_l)}}, \ \text{if} \ 1\le i,j\le k\\ 
&&\nonumber (P_k)_{ii}:=1, \ \ \text{if} \ i>k, \\ 
&&\nonumber (P_k)_{ij}:=0, \ \ \text{otherwise},
\end{eqnarray}
and its inverse is,
\begin{eqnarray}\label{clsQ}
&&(P_k^{-1})_{ij}=\frac{{(-1)^{i+j}}\Delta^{1,...,\hat{j},...,k}_{1,...,k-1}\left(A-\lambda^{(k)}_i\right)}{\sqrt{\prod_{l=1,l\ne i}^k(\lambda^{(k)}_i-\lambda^{(k)}_l)\prod_{l=1}^{k-1}(\lambda^{(k)}_i-\lambda^{(k-1)}_l)}}, \ \text{if} \ 1\le i,j\le k\\ 
&&\nonumber (P_k^{-1})_{ii}=1, \ \ \text{if} \ i>k, \\ 
&&\nonumber (P_k^{-1})_{ij}=0, \ \ \text{otherwise}.
\end{eqnarray}
Therefore, by definition the function
\begin{equation}\label{akiexp}
    a^{(k)}_i(A)= \sum_{v=1}^n(P_k(A)^{-1})_{iv} \cdot (A)_{v,k+1}=\frac{(-1)^{k+i}\Delta^{1,...,k}_{1,...,k-1,k+1}\left(A-\lambda^{(k)}_i\right)}{\sqrt{\prod_{l=1,l\ne i}^k(\lambda^{(k)}_i-\lambda^{(k)}_l)\prod_{l=1}^{k-1}(\lambda^{(k)}_i-\lambda^{(k-1)}_l)}}.
\end{equation}
Let $L^{(k+1)}(A)\in {\rm U}(k+1)\subset {\rm U}(n)$ be the matrix given by \begin{eqnarray}\label{L}
&&L^{(k+1)}_{ij}(A):=\frac{a^{(k)}_i}{ N_j^{(k+1)}(\lambda^{(k)}_i-\lambda^{(k+1)}_j)}, \ \text{for} \ i\ne k+1, \ j=1,...,k+1,\\ \label{L1}
&&L^{(k+1)}_{k+1,j}(A):=\frac{1}{N_j^{(k+1)}}, \ \ \ \ \ \ \text{for} \ j=1,...,k+1,
\end{eqnarray}
where the normalizer \begin{equation}\label{normalizer}N_j^{(k+1)}(A):=\sqrt{1+\sum_{l=1}^{k}\frac{|a_l^{(k)}|^2}{(\lambda^{(k)}_l-\lambda^{(k+1)}_j)^2}}.
\end{equation}
The upper left $k+1$-th submatrix of $L^{(k+1)}(A)^{-1}A_{k}L^{(k+1)}(A)$ is ${\rm diag}(\lambda^{(k+1)}_1,...,\lambda^{(k+1)}_{k+1})$ and the entries in $k+1$-th row of $L^{(k+1)}(A)$ is real and positive, thus we also have
\begin{eqnarray}\label{diagP}
P_{k+1}(A)=P_k(A)\cdot L^{(k+1)}(A).
\end{eqnarray}

Computing the character polynomial of the upper left $k+1$-th submatrix of $A_k$ leads to the identity 
\begin{equation}\label{normalizer1}
    |a_i^{(k)}|^2=-\frac{\prod_{v=1}^{k+1}(\lambda^{(k)}_i-\lambda^{(k+1)}_v)}{\prod_{v\ne i}^{k}(\lambda^{(k)}_i-\lambda^{(k)}_v)}.
\end{equation}
Thus the normalizer $N_j^{(k+1)}(A)$ only depends on the eigenvalues ${\lambda^{(k)}_j}$'s and ${\lambda^{(k+1)}_j}$'s, and can be rewritten as (using a combinatorial identity)
\begin{equation}\label{normalizer2}
N_j^{(k+1)}=\sqrt{\frac{\prod_{v=1, v\ne j}^{k+1}(\lambda^{(k+1)}_j-\lambda^{(k+1)}_v) }{\prod_{v=1}^{k}(\lambda^{(k+1)}_j-\lambda^{(k)}_v)}}
\end{equation}

\subsection{Riemann-Hilbert-Birkhoff maps at caterpillar points are compatible with Gelfand-Tsetlin systems}\label{rRHvsGZ}
Recall that in Section \ref{GTcoor}, we have introduced the Gelfand-Tsetlin maps and Thimm torus actions. In this subsection, we first introduce their multiplicative analogs \cite{FR} on $\Herm^+(n)$, and then prove that the map $\nu(u_{\rm cat})$ is compatible with them.

{\bf Logarithmic Gelfand-Tsetlin maps.} 
Let $\Herm^+(n)\subset \Herm(n)$ denote the subset of positive definite
Hermitian matrices, and define a logarithmic Gelfand-Tsetlin map 
\begin{equation}\label{eq:logmomentmap}
\mu\colon \Herm^+(n)\to \mathbb{R}^{\frac{n(n+1)}{2}},
\end{equation}
taking $A$ to the collection of numbers
$\mu^{(k)}_i(A)=\log(\lambda^{(k)}_i(A))$. Here recall that $\lambda^{(k)}_i(A)'s$ are the ordered set of
eigenvalues of the upper left $k$-th submatrix $A^{(k)}$ of $A\in \Herm(n)$. Then $\mu$ is a continuous
map from $\Herm^+(n)$ onto the Gelfand-Tsetlin cone $\mathcal{C}(n)$.

{\bf Thimm torus actions.}
Let $\mathcal{C}_0(n)\subset \mathcal{C}(n)$ denote the subset where all of the eigenvalue inequalities
\eqref{eq:cone} are strict. Let $\Herm_0^+(n)$ denote the intersection of $\Herm_0(n)$ and $\Herm^+(n)$, i.e., \[\Herm_0^+(n)=\mu^{-1}(\mathcal{C}_0(n)).\] Then the actions of Thimm torus on 
$\Herm_0(n)$ defined in \eqref{eq:taction} restrict to a torus action on $\Herm_0^+(n)$. The action preserves the logarithmic Gelfand-Tsetlin map $\mu$. 

Recall from Definition \ref{RHBcat} that the Riemann-Hilbert map at $u_{\rm cat}$ is
\[\nu(u_{\rm cat})\colon\Herm(n)\rightarrow \Herm^+(n); \ A\mapsto S_-(u_{\rm cat},A) S_+(u_{\rm cat},A).\]
\begin{pro}\cite[Proposition 4.1]{Xu}\label{ConnGZ}
The Riemann-Hilbert map 
\[\nu(u_{\rm cat}): {\rm Herm}(n)\rightarrow {\rm Herm}^+(n)\] 
is a diffeomorphism compatible with the Gelfand-Tsetlin systems. That is
\begin{itemize}
\item[(a).] $\nu(u_{\rm cat})$ intertwines the
  Gelfand-Tsetlin maps: $\mu\circ \nu(u_{\rm cat})=\lambda$;
\item[(b).] $\nu(u_{\rm cat})$ intertwines the Thimm's torus actions on $\Herm_0(n)$ and $\Herm_0^+(n)$. 
\end{itemize}
\end{pro}

\subsection{Explicit formula of the Alekseev-Meinrenken diffeomorphism}\label{AMdiff}
In \cite{AM}, the authors introduced a distinguished diffeomorphism $\Gamma_{AM}:\Herm(n)\rightarrow\Herm^+(n)$ compatible with the Gelfand-Tsetlin systems. In this subsection, we point out the relation between $\Gamma_{AM}$ and $\nu(u_{\rm cat})$, and use the expression of Stokes matrices at a caterpillar point to derive an explicit expression of $\Gamma_{AM}$.

\subsubsection{The Alekseev-Meinrenken diffeomorphim}
The Poisson manifolds ${\frak u}(n)^*\cong\Herm(n)$ and ${\rm U}(n)^*\cong\Herm^+(n)$ carry the structures of Gelfand–Tsetlin integrable system \cite{GS} and multiplicative Gelfand-Tsetlin system \cite{FR} respectively. The transpose map $T:A\mapsto A^T$ is an anti-Poisson
involution of $\Herm(n)$, as well as of $\Herm^+(n)$. The fixed point set of the transpose map on $\Herm(n)$ is the set $\Sym(n)$ of real symmetric $n$ by $n$ matrices, and the intersection $\Sym_0(n)=\Sym(n)\cup \Herm_0(n)$ is a union of Lagrangian submanifolds of $\Herm_0(n)$. Actually the space $\Sym_0(n)$ has $2^{\frac{n(n-1)}{2}}$ components.

\begin{thm}\cite{AM}\label{111}
There exists a unique diffeomorphism
$\Gamma_{AM}: \Herm(n)\rightarrow \Herm^+(n)$
with the following three properties:

$(a)$ $\Gamma_{AM}$ intertwines the Gelfand-Tsetlin maps: $\mu\circ\Gamma_{AM}=\lambda$;

$(b)$ $\Gamma_{AM}$ intertwines the Thimm's torus actions on $\Herm_0(n)$ and $\Herm^+_0(n)$;

$(c)$ For any connected component $S$ of $\Sym_0(n)\subset \Herm(n)$, $\Gamma_{AM}(S)\subset S$.

In particular, the map $\Gamma_{AM}$ is a Poisson isomorphism.
\end{thm}

\subsubsection{Principal bundles and cross sections}\label{principalbundle}
\begin{pro}(see e.g., \cite{AM} Proposition 2.1)\label{prinbundle}
The restriction of the Gelfand-Tsetlin map $\lambda$ to $\Herm_0(n)$ defines a principal bundle
$$\lambda: \Herm_0(n) \rightarrow \mathcal{C}_0(n)$$
over the cone $\mathcal{C}_0(n)$ with structure group the Thimm's torus $T^{\frac{n(n-1)}{2}}$. 
It further restricts
to a principal bundle
$$\lambda: \Sym_0(n)\rightarrow \mathcal{C}_0(n)$$
with discrete structure group $T_\mathbb{R}^{\frac{n(n-1)}{2}}\cong (\mathbb{Z}_2)^{\frac{(n-1)n}{2}}$. 
Similarly for (the restriction of) the logarithmic
Gelfand-Tsetlin map $\mu: \Herm_0^+(n)\rightarrow \mathcal{C}_0(n)$ and $\mu:\Sym_0^+(n)\rightarrow \mathcal{C}_0(n)$. Here $\Sym_0^+(n)$ denotes the intersection of $\Herm_0^+(n)$ with $\Sym(n)$.
\end{pro}
Note that $\Sym_0(n)$ has $2^{\frac{n(n-1)}{2}}$ components, and the "real part" $T_\mathbb{R}^{\frac{n(n-1)}{2}}\cong (\mathbb{Z}_2)^{\frac{(n-1)n}{2}}$ of the Thimm's torus action on $\Sym_+(n)$ relates the different connected components. 
Any connected component $S$ of $\Sym_0(n)\subset \Herm_0(n)$ can be understood as a cross section of the principal bundle $\lambda: \Herm_0(n) \rightarrow \mathcal{C}_0(n)$. Therefore, a principal $T^{\frac{n(n-1)}{2}}$-bundle map $$\Gamma:(\Herm_0(n),\lambda)\rightarrow (\Herm_0^+(n),\mu={\rm log}\circ \lambda)$$
is uniquely characterized by the image $\Gamma(S)\subset \Herm^+_0(n)$ of $S$. In particular, we can specify a bundle map $\Gamma_0$ by imposing the real condition: $\Gamma_0(S)\subset S$ for a given connected component $S$ of $\Sym_0(n)$. By the $T_\mathbb{R}^{n(n-1)/2}$-equivariance, $\Gamma_0(S)\subset S$ holds true for all connected component $S$ of $\Sym_0(n)$. Thus it gives the geometric interpretation of the conditions $(a)$-$(c)$ characterising the Alekseev-Meinrenken diffeomorphism $\Gamma_{AM}$ in Theorem \ref{11}. Actually, the map $\Gamma_0$ coincides with the restriction of $\Gamma_{AM}$ on $\Herm_0(n)$. However, the existence of a smooth extension of $\Gamma_0$ to $\Herm(n)$ is not obvious.

Therefore, the maps $\Gamma_{AM}$ and $\nu(u_{\rm cat})$ are both Poisson diffeomorphism intertwining the Gelfand-Tsetlin systems, while map a Lagrangian submanifold $S\subset \Sym_0(n)$ to two different cross sections of the Gelfand-Tsetlin map $\mu$.

\subsubsection{Phase transformations}\label{Phasetrans}
Let $\theta:\Herm_0(n)\rightarrow T^1\times \cdot\cdot\cdot\times T^{n-1}=T^{\frac{n(n-1)}{2}}$ be a map from $\Herm_0(n)$ to the Gelfand-Tsetlin tours. Associated to $\theta$, there is a gauge transformation \begin{equation}\label{gauge}
\mathcal{X}_{\theta}:\Herm_0(n)\rightarrow \Herm_0(n); \ A\mapsto {\theta(A)}\bullet A,
\end{equation}
where $\theta(A)\bullet A$ denotes the Thimm's action of $\theta(A)$ on $A$, see \eqref{eq:taction}.
For each $1\le k\le n-1$, let us denote the $T^k$-component of $\theta$ by
\[\theta^{(k)}(A)={\rm diag}(e^{\I\theta_1^{(k)}(A)},...,e^{\I\theta_{k}^{(k)}(A)})\in T^{k}, \ \ \text{ for any } \ A\in\Herm_0(n).\]
Under the Gelfand-Tsetlin action and angle coordinates, the Thimm action of the element \[\theta^{(k)}(A)={\rm diag}(e^{\I \theta^{(k)}_1},...,e^{\I \theta^{(k)}_k})\in T(k)\] on a point $p$ with coordinates $(\{\lambda^{(i)}_j\}_{1\le j\le i\le n}, \{\psi^{(i)}_j\}_{1\le j\le i\le n-1} )$ is described by
\begin{eqnarray}\label{torusGT}
\lambda^{(i)}_j(\theta^{(k)}(A)\bullet p) =\lambda^{(i)}_j(A), \hspace{5mm} \psi^{(i)}_j(\theta^{(k)}(A)\bullet p) =\psi^{(i)}_j(A)+\delta_{ki}\theta^{(k)}_{j}(A).
\end{eqnarray}
Then it follows from \eqref{torusGT} that the Gelfand-Tsetlin action and angle coordinates on $\Herm_0(n)$ change under the map $\mathcal{X}_{\theta}$ as \begin{eqnarray}\label{action}
\lambda^{(i)}_j(\mathcal{X}_{\theta}(A))=\lambda^{(i)}_j(A), \hspace{5mm} \psi^{(i)}_j(\mathcal{X}_{\theta}(A))=\psi^{(i)}_j(A)+\sum_{k=1}^{n-1}\delta_{ik}\theta_j^{(k)}(A).
\end{eqnarray}

\subsubsection{Proof of Theorem \ref{exAM}}\label{sec:exAM}
We have seen that the maps $\Gamma_{AM}$ and $\nu(u_{\rm cat})$ relate to each other by a gauge transformation. Thus to prove Theorem \ref{exAM}, we only need to construct explicitly the transformation. We decompose it into three steps.

\vspace{2mm}
{\bf Step 1. Reformulation of the Riemann-Hilbert map $\nu(u_{\rm cat})$ at the caterpillar point $u_{\rm cat}$.} 
Recall that $L^{(j)}(A)$ is defined in \eqref{L} for all $1\le j\le n-1$.
For any integer $1\le k\le n$ and any $A\in\Herm_0(n)$, we introduce the {\it normalized connection matrix}
\begin{eqnarray}\label{tildeCk}
\widetilde{C}\left(E_k,\delta_k(A_{{k-1}})\right):=C\left(E_k, \delta_k(A_{k-1})\right)\cdot L^{(k)}(A).
\end{eqnarray}
Here recall from \eqref{zhuazi} that $A_{{k-1}}=P_{k-1}(A)^{-1}AP_{k-1}(A)$ is the diagonalization of the upper left $(k-1)$-th submatrix of $A$.
\begin{lem}\label{mono2}
For any $A\in {\rm Herm}_0(n)$, if we define \begin{equation}\label{tildeC}
    \widetilde{C}(u_{\rm cat},A)=\widetilde{C}(E_2,\delta_2(A_1)) \widetilde{C}(E_3,\delta_3(A_2))\cdot\cdot\cdot \widetilde{C}(E_n,\delta_n(A)) \in U(n),
\end{equation}
as the pointwise multiplication, then we have the Riemann-Hilbert map at $u_{\rm cat}$ (with respect to $U_{\rm id}$)
\begin{equation}\label{productsC}
\nu(u_{\rm cat},A)=\widetilde{C}(u_{\rm cat}, A)e^{A_{n}}\widetilde{C}(u_{\rm cat},A)^{-1}.
\end{equation}
Here recall that $A_{n}={\rm diag}(\lambda_1^{(n)},...,\lambda^{(n)}_n)$. 
\end{lem}
\pf 
For any $1\le k\le n$, we take the obvious inclusion of $U(k)$ as the upper left corner of $U(n)$, extended
by $1$'s along the diagonal. Since $U(k-1)\in U(n)$ is in the centralizer of the irregular term $\I  E_k$ of the equation \begin{eqnarray}\label{eq:rel}
\frac{dF}{dz}=\left(\I  E_k-\frac{1}{2\pi \I }\frac{\delta_k(A)}{z}\right)F,
\end{eqnarray}
the connection matrix $C(E_k,\delta_k(A))$ of \eqref{eq:rel} has the following ${U}(k-1)$-equivariance:
for any $G\in {U}(k-1)\subset {U}(n)$, 
\[C\left(E_k, G\delta_k(A)G^{-1}\right)=G\cdot C\left(E_k,\delta_k(A)\right)\cdot G^{-1}.\]
It follows from the property of equivariance, the identity \eqref{diagP} and Definition \eqref{tildeC} that \[\widetilde{C}(u_{\rm cat},A)=\left(\overrightarrow{\prod_{k=1,...,n}}C\left(E_k,\delta_k(A)\right)\right)\cdot P_n(A).\] Then the proposition follows from the identity \eqref{nuC}, and the identity $P_n(A)e^{A_{n}}P_n(A)^{-1}=e^A$ (recall the definition $A_{n}=P_n(A)^{-1}AP_n(A)$ of $A_n$).
\qed

\vspace{2mm}
Based on \eqref{productsC} and \eqref{tildeC}, to get the explicit expression of $\nu(u_{\rm cat},A)$, one only needs to compute the normalized connection matrix $\widetilde{C}\left(E_{k+1},\delta_{k+1}(A_{k})\right)$ given in \eqref{tildeCk}. In particular, the following proposition can be found in \cite{Balser}, see also \cite{Xu1}.
\begin{pro}\label{explicitmCS}
The entries of the normalized connection matrix $\widetilde{C}\left(E_{k+1},\delta_{k+1}(A_{k})\right)$ are given by\begin{eqnarray*} \widetilde{C}_{ij}=\frac{ e^{\frac{\lambda^{(k)}_i-\lambda^{(k+1)}_j}{4}}}{(\lambda^{(k)}_i-\lambda^{(k+1)}_j)}\frac{\prod_{v=1}^{k+1}\Gamma\left(1+\frac{\lambda^{(k+1)}_v-\lambda^{(k+1)}_j}{2\pi \I }\right)\prod_{v=1}^{k}\Gamma\left(1+\frac{\lambda^{(k)}_v-\lambda^{(k)}_i}{2\pi \I }\right)}{\prod_{v=1,v\ne i}^{k}\Gamma\left(1+\frac{\lambda^{(k)}_v-\lambda^{(k+1)}_j}{2\pi \I }\right)\prod_{v=1,v\ne j}^{k+1}\Gamma\left(1+\frac{\lambda^{(k+1)}_v-\lambda^{(k)}_i}{2\pi \I }\right)}\cdot \frac{a^{(k)}_i(A)}{ N_j^{(k+1)}(A)},
\end{eqnarray*}
for $1\le j\le k+1, 1\le i\le k$, and
\begin{eqnarray*} \widetilde{C}_{k+1,j}&=&\frac{e^{\frac{\lambda^{(k+1)}_j-(A)_{k+1,k+1}}{4}}\prod_{v=1}^{k+1}\Gamma\left(1+\frac{\lambda^{(k+1)}_v-\lambda^{(k+1)}_j}{2\pi \I }\right)}{N_j^{(k+1)}(A)\cdot \prod_{v=1}^{k}\Gamma\left(1+\frac{\lambda^{(k)}_v-\lambda^{(k+1)}_j}{2\pi \I }\right)}, \ \ \ \text{for} \ 1\le j\le k+1.\\
\widetilde{C}_{ii}&=&1, \ \ \ \ \ \ \text{for} \ k+1<i\le n,\\
\widetilde{C}_{ij}&=&0, \ \ \ \ \ \ \text{otherwise}.
\end{eqnarray*}
\end{pro}

From the expression in Proposition \ref{explicitmCS}, we observe that for any $1<i\le n$ the normalized connection matrix has a decomposition
\begin{equation}\label{decomC}
\widetilde{C}\left(E_{i},\delta_{i}(A)\right)=D^{(i-1)}_L(A)\cdot {\rm diag}(-a^{(i-1)}_1,...,-a^{(i-1)}_{i-1},1,...,1)\cdot R^{(i)}(A)\cdot D^{(i)}_R(A),
\end{equation}
where $D^{(i-1)}_L(A)$ is the $n$ by $n$ diagonal matrix
\begin{eqnarray*}
D^{(i-1)}_{L,kk}(A)=\frac{\prod_{v=1}^{i-1}\Gamma\left(1+\frac{\lambda^{(i-1)}_v-\lambda^{(i-1)}_k}{2\pi\I }\right)}{\prod_{v=1}^{i}\Gamma\left(1+\frac{\lambda^{(i)}_v-\lambda^{(i-1)}_k}{2\pi\I }\right)}, \ \ \text{for} \ 1\le k\le i-1,\ \ \ \text{and} \ \ D^{(i-1)}_{L,ii}(A)=1,
\end{eqnarray*}
$D^{(i)}_R(A)$ is the $n$ by $n$ diagonal matrix function
\begin{eqnarray*}
D^{(i)}_{R,kk}(A)=\frac{ \prod_{v=1}^{i}\Gamma\left(1+\frac{\lambda^{(i)}_v-\lambda^{(i)}_k}{2\pi\I }\right)}{\prod_{v=1}^{i-1}\Gamma\left(1+\frac{\lambda^{(i-1)}_v-\lambda^{(i)}_k}{2\pi\I }\right)}, \ \text{for} \ 1\le k\le i,
\end{eqnarray*}
and $R^{(i)}(A)$ is a $n$ by $n$ matrix with entries

\begin{equation}\label{RA}
R^{(i)}_{kj}(A)=\left\{
          \begin{array}{lr}
            \frac{e^{\frac{\lambda^{(i-1)}_k-\lambda^{(i)}_j}{4}}}{N_j^{(i)}{\rm sinh}\Big(\frac{\lambda^{(i-1)}_k-\lambda^{(i)}_j}{2}\Big)},   & \text{if} \ \ 1\le k,j\le i\\
           1, & \text{if} \ \  k=j\in\{i+1,i+2,...,n\} \\
           0, & \ \text{otherwise},
             \end{array}
\right. 
\end{equation}

Now let us denote by 
\begin{eqnarray}\label{Diag}
D^{(i)}(A):=D^{(i)}_R(A) \cdot D^{(i)}_L(A)\cdot {\rm diag}(-a^{(i)}_1,...,-a^{(i)}_{i},1,...,1)
\end{eqnarray} the multiplication of the three diagonal matrices. Then we can rewrite the map $\widetilde{C}(u_{\rm cat},A)$ in Lemma \ref{mono2} as \[\left(\overrightarrow{\prod_{i=1,...,n-1}}D^{(i)}(A)R^{(i)}(A)\right)\cdot D^{(n)}_R(A).\] 
Since $D^{(n)}_R(A)$ commutes with $e^{A_{n}}$ (recall that $A_{n}={\rm diag}(\lambda^{(n)}_1,...,\lambda^{(n)}_n)$), we see that 
\begin{cor}\label{rBformula}
The Riemann-Hilbert-Birkhoff map $\nu(u_{\rm cat})$ can be rewritten as
\begin{eqnarray}
\nu(u_{\rm cat}):{\rm Herm}_0(n)\rightarrow {\rm Herm}^+_0(n); \ A\mapsto \mathcal{T}(A)e^{A_{n}}\mathcal{T}(A)^{-1}
\end{eqnarray}
where $\mathcal{T}(A):=D^{(1)}(A)R^{(1)}(A)\cdot\cdot\cdot D^{(n-1)}(A)R^{(n-1)}(A)$. 
\end{cor}

{\bf Step 2. The phase transformation of $\nu(u_{\rm cat})$.} Based on the expression \eqref{action} and Corollary \ref{rBformula}, we have
\begin{pro}\label{RealrB}
There exists a unique map $\theta:\Herm_0(n)\rightarrow T^{\frac{n(n-1)}{2}}$, such that for all $A\in\Sym_0$ with the angle coordinates $\psi^{(i)}_j(A)=0$, the Gelfand-Tsetlin angle variables 
\begin{equation}\label{psireal}
\psi^{(i)}_j\left(\nu(u_{\rm cat}, \mathcal{X}_{\theta}(A)\right)=0, \ \ \ \text{for \ all} \ 1\le j\le i\le n-1.
\end{equation}
\end{pro}
\pf It follows from Corollary \ref{rBformula} that for any $\theta$, the composed map
\begin{eqnarray}\label{phasechange}
\nu(u_{\rm cat})\circ \mathcal{X}_{\theta}(A):{\rm Herm}_0(n)\rightarrow {\rm Herm}^+    _0(n); \ \mathcal{T}(\mathcal{X}_{\theta}(A))e^{A_{d_n}}\mathcal{T}(\mathcal{X}_{\theta}(A))^{-1}
\end{eqnarray}
where $\mathcal{T}(\mathcal{X}_{\theta}(A))=D^{(1)}(\mathcal{X}_{\theta}(A))R^{(1)}(\mathcal{X}_{\theta}(A))\cdot\cdot\cdot D^{(n-1)}(\mathcal{X}_{\theta}(A))R^{(n-1)}(\mathcal{X}_{\theta}(A))$. 

Now let us compare the difference between $\mathcal{T}(\mathcal{X}_{\theta}(A))$ and $\mathcal{T}(A)$ using the structure of the decomposition \eqref{decomC}.
Following \eqref{Diag}, \[D^{(i)}(A)=D^{(i)}_R(A) \cdot D^{(i)}_L(A)\cdot {\rm diag}\left(-a^{(i)}_1(A),...,-a^{(i)}_{i}(A),1,...,1\right)\] for any $1\le i\le n-1$. By the expression \eqref{RA}, $R^{(i)}(A)$ only depends on the Gelfand-Tsetlin action variables, 
which are preserved under the action $\mathcal{X}_{\theta}$. Thus $R^{(i)}(A)=R^{(i)}(\mathcal{X}_{\theta}(A))$. 
For the same reason, we have $D^{(i)}_L(A)=D^{(i)}_L(\mathcal{X}_{\theta}(A))$ and $D^{(i)}_R(A)=D^{(i)}_R(\mathcal{X}_{\theta}(A))$. Therefore, 
\[\mathcal{T}(\mathcal{X}_{\theta}(A))=\overrightarrow{\prod_{i=1,...,n-1}}\left(D^{(i)}_R(A) \cdot D^{(i)}_L(A)\cdot {\rm diag}\left(-a^{(i)}_1(\mathcal{X}_{\theta}(A)),...,-a^{(i)}_{i}(\mathcal{X}_{\theta}(A)),1,...,1\right) R^{(i)}(A)\right),\]
which only differs from $\mathcal{T}(A)$ by the angle of the variables $a^{(i)}_j$. For the angles $\psi^{(i)}_j$ of the variables $a^{(i)}_j$, by \eqref{action} we have $\psi^{(i)}_j(\mathcal{X}_{\theta}(A))=\theta_j^{(i)}(A)+\psi^{(i)}_j(A)+\pi \I $.
It implies
\begin{eqnarray*}
D^{(i)}(\mathcal{X}_{\theta}(A))=D^{(i)}_R(A) D^{(i)}_L(A)\cdot {\rm diag}\left(e^{\theta_1^{(i)}(A)+\psi^{(i)}_1(A)+\pi\I },...,e^{\theta_i^{(i)}(A)+ \psi^{(i)}_{i}(A)+\pi\I },1,...,1\right).
\end{eqnarray*}
The $(j,j)$ entry of the above diagonal matrix is 
\begin{eqnarray*}
&&{\rm Arg}(D^{(i)}_{jj}(\mathcal{X}_{\theta}(A)))={\rm Arg}(D^{(i)}_{R,jj}(A)  D^{(i)}_{L,jj}(A))+\psi^{(i)}_j(A)+\pi\I +\theta_j^{(i)}(A), \ j\le i\\
&&{\rm Arg}(D^{(i)}_{jj}(\mathcal{X}_{\theta}(A)))=1, \ j>i.
\end{eqnarray*}

The above discussion is for any $\theta$. Let us now choose the particular $\theta$ such that \eqref{psireal} holds. For any $A\in\Sym_0$ with the angle coordinates $\psi^{(i)}_j(A)=0$, the condition ${\rm Arg}(D^{(i)}_{jj}(\mathcal{X}_{\theta}(A)))=0$, for all $1\le i\le n-1$ and $1\le j\le i$, is equivalent to a linear system 
\begin{equation}\label{linearsys}
    \theta_j^{(i)}(A)=-{\rm Arg}(D^{(i)}_{R,jj}(A) D^{(i)}_{L,jj}(A))-\pi\I , \ \ \text{for all} \ 1\le j\le i\le n-1.
\end{equation}
The linear system has a unique solution, i.e., a collection of $\theta_j^{(i)}(A)\in [0,2\pi)$. 
Let us denote the solution of the linear system by $\theta:\Herm_0(n)\rightarrow T^{\frac{n(n-1)}{2}}$. That is $\theta(A)$ is the collection of $\theta_j^{(i)}(A)$ satisfying \eqref{linearsys} for all $1\le j\le i\le n-1$.

From the above discussion, for any $A\in\Sym_0$ the matrix $D^{(i)}(\mathcal{X}_{\theta}(A))R^{(i)}(\mathcal{X}_{\theta}(A))$ is real. Together with \eqref{phasechange}, we see that the composed map $\nu(u_{\rm cat})\circ \mathcal{X}_{\theta}$ restricts to a map from $\Sym_0$ to $\Sym_0^+$. If we impose further the angle coordinates $\phi^{(i)}_j(A)=0$, then the angle variables $\phi^{(i)}_j$ of $\nu(u_{\rm cat})\circ \mathcal{X}_{\theta}(A)$ is also zero. It finishes the proof. \qed

\vspace{2mm}
{\bf Step 3. The explicit expression of $\Gamma_{AM}$.} Now it is straightforward to write down the map $\Gamma_{AM}$. First by the definition of $\Gamma_{AM}$ and Proposition \ref{RealrB}, we have 
\begin{cor}\label{AMcat}
The map $\mathcal{X}_\theta$ extends to a diffeomorphism from $\Herm_0(n)$ to $\Herm(n)$ and is such that $\Gamma_{AM}=\nu(u_{\rm cat})\circ \mathcal{X}_\theta$. That is
\[\Gamma_{AM}(A)=\nu(u_{\rm cat},\mathcal{X}_\theta(A))=\widetilde{C}(u_{\rm cat},\mathcal{X}_\theta(A))\cdot A_n\cdot \widetilde{C}(u_{\rm cat},\mathcal{X}_\theta(A))^{-1}, \ \ \forall \ A\in\Herm(n).\]
\end{cor}
Following the proof of Proposition \ref{RealrB}, the phase transformation $\theta$ transforms the gamma functions of the form $\Gamma(1+\frac{r}{2\pi\I })$ for $r\in\mathbb{R}$, in the expression of $\widetilde{C}(u_{\rm cat},\mathcal{X}_\theta(A))$ for any $A\in\Sym_0$, to their real parts. By the complex conjugate of gamma function and Euler's reflection formula
\begin{eqnarray*}
&&\overline{\Gamma(z)}=\Gamma(\bar{z}), \ z\notin \mathbb{C}\setminus\{0,-1,-2,...\}, \\ 
&&\Gamma(z)\Gamma(1-z)=\frac{\pi}{{\rm sin}(\pi z)}, \hspace{3mm} z\notin \mathbb{Z},
\end{eqnarray*}
we get $|\Gamma(1+\frac{r}{2\pi\I })|=\sqrt{\frac{2r }{2{\rm sinh}(\frac{r}{2})}}$ for $r\in\mathbb{R}$. Thus replacing the terms $\Gamma(1+\frac{r}{2\pi\I })$ by their norms $\sqrt{\frac{2r}{2{\rm sinh}(\frac{r}{2})}}$, we get the expression of the matrix $\widetilde{C}\left(E_k,\delta_k(\mathcal{X}_\theta(A))\right)$,
\begin{eqnarray*}
\widetilde{C}\left(E_k,\delta_k(\mathcal{X}_\theta(A))\right)_{ij} = \frac{e^{\frac{\lambda^{(k-1)}_i-\lambda^{(k)}_j}{4}}a^{(k-1)}_i}{N_j^{(k)}(\lambda^{(k-1)}_i-\lambda^{(k)}_j)} \sqrt{\frac{H_{ij}^{(k)} \prod_{v=1,v\ne i}^{k-1}{\rm sinh}\Big(\frac{\lambda^{(k-1)}_v-\lambda^{(k)}_j}{2}\Big)\prod_{v=1,v\ne j}^{k}{\rm sinh}\Big(\frac{\lambda^{(k-1)}_i-\lambda^{(k)}_v}{2}\Big)}{\prod_{v=1, v\ne j}^{k}{\rm sinh}\Big(\frac{\lambda^{(k)}_v-\lambda^{(k)}_j}{2}\Big)\prod_{v=1, v\ne i}^{k-1}{\rm sinh}\Big(\frac{\lambda^{(k-1)}_i-\lambda^{(k-1)}_v}{2}\Big)}}, 
\end{eqnarray*}
for $1\le i\le k-1, 1\le j\le k$, and
\begin{eqnarray*} && \widetilde{C}\left(E_k,\delta_k(\mathcal{X}_\theta(A))\right)_{kj}=\frac{e^{\frac{\lambda^{(k)}_j-\lambda^{(k-1)}_k}{4}}}{N_j^{(k)}}\sqrt{\frac{\prod_{v=1}^{k-1}{\rm sinh}\Big(\frac{\lambda^{(k-1)}_v-\lambda^{(k)}_j}{2}\Big)\prod_{v=1,v\ne j}^{k}(\lambda^{(k)}_v-\lambda^{(k)}_j)}{\prod_{v=1, v\ne j}^{k}{\rm sinh}\Big(\frac{\lambda^{(k)}_v-\lambda^{(k)}_j}{2}\Big)\prod_{v=1}^{k-1}(\lambda^{(k-1)}_v-\lambda^{(k)}_j)}}, \ \text{for}  \ 1\le j\le k,\\
&&\widetilde{C}\left(E_k,\delta_k(\mathcal{X}_\theta(A))\right)_{ii}=1, \ \ \ \ \ \ \text{for} \ k<i\le n,\\
&&\widetilde{C}\left(E_k,\delta_k(\mathcal{X}_\theta(A))\right)_{ij}=0, \ \ \ \ \ \ \text{otherwise}.
\end{eqnarray*}
Here $N_j^{(k)}$ and $H^{(k)}_{ij}$ are given by \begin{eqnarray*}
N_j^{(k)}:=\sqrt{1+\sum_{l=1}^{k-1}\frac{|a_l^{(k-1)}|^2}{(\lambda^{(k-1)}_l-\lambda^{(k)}_j)^2}}, \ \ \ \ \ \ \label{H}H^{(k)}_{ij}=\frac{\prod_{v=1, v\ne j}^{k}(\lambda^{(k)}_v-\lambda^{(k)}_j)\prod_{v=1, v\ne i}^{k-1}(\lambda^{(k-1)}_i-\lambda^{(k-1)}_v)}{\prod_{v=1,v\ne i}^{k-1}(\lambda^{(k-1)}_v-\lambda^{(k)}_j)\prod_{v=1,v\ne j}^{k}(\lambda^{(k-1)}_i-\lambda^{(k)}_v)}.
\end{eqnarray*}
Using the identities \eqref{normalizer1} and \eqref{normalizer2}, we have \[|a_i^{(k-1)}|^2=-\frac{\prod_{v=1}^{k}(\lambda^{(k-1)}_i-\lambda^{(k)}_v)}{\prod_{v\ne i}^{k-1}(\lambda^{(k-1)}_i-\lambda^{(k-1)}_v)}, \ \ N_j^{(k)}=\sqrt{\frac{\prod_{v=1, v\ne j}^{k}(\lambda^{(k)}_j-\lambda^{(k)}_v) }{\prod_{v=1}^{k-1}(\lambda^{(k)}_j-\lambda^{(k-1)}_v)}}.\] 
In the end, let us check that $\widetilde{C}\left(E_k,\delta_k(\mathcal{X}_\theta(A))\right)_{ij}$ coincides with the expression of $\psi^{(k)}(A)_{ij}$ in Theorem \ref{exAM}. To see this, first notice that the identity $\sqrt{H_{ij}^{(k)}}=N_j^{(k)}(\lambda^{(k-1)}_i-\lambda^{(k)}_j)/|a^{(k-1)}_i|$ (the identity itself follows from the character polynomial of the upper left $k$-th submatrix of $A_k$). Then the identities \eqref{normalizer1} and \eqref{akiexp} imply 
\[\frac{a^{(k-1)}_i\sqrt{H_{ij}^{(k)}}}{N_j^{(k)}(\lambda^{(k-1)}_i-\lambda^{(k)}_j)} =\frac{a^{(k-1)}_i}{|a^{(k-1)}_i|}=\frac{(-1)^{k-1+i}\Delta^{1,...,k-1}_{1,...,k-2,k}\left(A-\lambda^{(k-1)}_i\right)}{\sqrt{-{\prod_{l=1}^{k}(\lambda^{(k-1)}_i-\lambda^{(k)}_l)\prod_{l=1}^{k-2}(\lambda^{(k-1)}_i-\lambda^{(k-2)}_l)}}}.\]
It finishes the proof of Theorem \ref{exAM}.

\subsubsection{Example: 2 by 2 cases.}
Let $A=\left(
  \begin{array}{cc}
     a & b \\
    \bar{b} & c 
  \end{array}\right)$ be a 2 by 2 Hermitian matrix. We will denote by $$\{\lambda^{(1)}_1:=a, \lambda^{(2)}_1, \lambda^{(2)}_2, \psi^{(2)}_1:=b/{|b|}\}$$ the corresponding Gelfand-Tsetlin coordinates. Here $\lambda^{(2)}_1, \lambda^{(2)}_2$ are the eigenvalues of $A$.
In this case, the formula in Theorem \ref{mainthm} gives us $$\Gamma_{AM}:\Herm(2)\rightarrow \Herm^+(2); \ A=\left(
  \begin{array}{cc}
     a & b \\
    \bar{b} & c 
  \end{array}\right)\mapsto \left(
  \begin{array}{cc}
     a' & b' \\
    \bar{b'} & c' 
  \end{array}\right),$$
where
$$
\left\{
  \begin{array}{ll}
a'=e^{\lambda^{(1)}_1}\\
b'=e^{i\psi^{(2)}_1}\sqrt{e^{\lambda^{(1)}_1+\lambda^{(2)}_1}+e^{\lambda^{(1)}_1+\lambda^{(2)}_2}-e^{2\lambda^{(1)}_1}-e^{\lambda^{(2)}_1+\lambda^{(2)}_2}}\\
c'=e^{\lambda^{(2)}_1}+e^{\lambda^{(2)}_2}-e^{\lambda^{(1)}_1}.
\end{array}
\right.
$$

The above expression in $n=2$ case coincide with the one given in \cite{AM}. However, for general $n$, coordinate expressions for the Alekseev-Meinrenken diffeomorphism or the Ginzburg-Weinstein maps were not known in the previous works \cite{Anton, AM, Boalch, EEM, GW}.

\subsection{Proof of Theorem \ref{mainthm}}\label{endsection}
On the one hand, replacing $A$ by $g(u;A)\cdot A\cdot g(u;A)^{-1}$ in Corollary \ref{AMcat}, we have
\[\Gamma_{AM}\left(\mathcal{X}_\theta^{-1}({\rm Ad}_{g(u;A)} A)\right)=\nu\left(u_{\rm cat},g(u;A)\cdot A\cdot g(u;A)^{-1}\right).\]
Here $\mathcal{X}_\theta^{-1}$ is the inverse map of $\mathcal{X}_\theta$. 
On the other hand, by Proposition \ref{leadingterm}, we have
\[\nu(u,A)=\nu\left(u_{\rm cat}, \hspace{2mm} g(u;A)\cdot A\cdot g(u;A)^{-1}\right)+\sum_{k=2}^{n-1}\mathcal{O}\left(\frac{u_{k}-u_{k-1}}{u_{k+1}-u_{k}}\right).\]
Combining the above two identities shows that if we introduce the diffeomorphism 
\[\psi=\mathcal{X}_\theta^{-1}\circ {\rm Ad}_{g(u;\cdot )}:A\mapsto \mathcal{X}_\theta^{-1}({\rm Ad}_{g(u;A)} A), \]
then $\psi$ is the diffeomorphism required in Theorem \ref{mainthm}, i.e., 
\[\nu(u,A)=\Gamma_{AM}(\psi(A))+\sum_{k=2}^{n-1}\mathcal{O}\left(\frac{u_{k}-u_{k-1}}{u_{k+1}-u_{k}}\right).\]

Now let us write down $\psi$ explicitly. First, for any $A$ and $u$, let us introduce an element in the product of torus
\begin{equation}\label{toriele}
(u_2-u_1)^{\frac{-\lambda^{(1)}(A)}{2\pi\I }}\times \overrightarrow{\underset{k=2,...,n-1}{\prod} }\left(\frac{u_{k+1}-u_{k}}{u_{k}-u_{k-1}}\right)^{\frac{-\lambda^{(k)}(A)}{2\pi\I }}\in T(1)\times \cdots \times T(n-1)\end{equation}
where $\lambda^{(k)}(A):={\rm diag}(\lambda^{(k)}_1,...,\lambda^{(k)}_{k})$. Let us introduce the diagonal matrix 
\begin{equation}\label{diagele}D(u;A)={\rm diag}\left(1,(u_2-u_1)^{-{A_{22}}/{2\pi\I }},...,(u_n-u_{n-1})^{-{A_{nn}}/{2\pi\I }}\right)\in T(n).
\end{equation}
Then one checks 
\begin{align*}
&g(u;A)Ag(u;A)^{-1}\\
=&D(u;A)\cdot \left((u_2-u_1)^{\frac{-\lambda^{(k)}(A)}{2\pi\I }}\times \overrightarrow{\underset{k=2,...,n-1}{\prod} }\left(\frac{u_{k+1}-u_{k}}{u_{k}-u_{k-1}}\right)^{\frac{-\lambda^{(k)}(A)}{2\pi\I }}\bullet A\right)\cdot D(u;A)^{-1}.
\end{align*}
Here $\bullet$ denotes the Thimm action. Therefore, under the Gelfand-Tsetlin action and angle coordinates in Definition \ref{GZfunctions} we have
\begin{align*}
\lambda^{(k)}_i({\rm Ad}_{g(u;A)}A)&=\lambda^{(k)}_i(A),\\
\psi^{(k)}_i({\rm Ad}_{g(u;A)}A)&=\psi^{(k)}_i(A)+ {\rm Arg}\left( (u_k-u_{k-1})^{\frac{\lambda^{(k)}_i(A)-A_{kk}}{2\pi\I }} (u_{k+1}-u_{k})^{\frac{A_{k+1,k+1}-\lambda^{(k)}_i(A)}{2\pi\I }}\right).
\end{align*}
Then the map ${\rm Ad}_{g(u;\cdot)}=\mathcal{X}_{\theta_1}$ (according to the definition in \eqref{gauge}) is a transformation along the Thimm torus fibration generated by a map \[\theta_1(u):\Herm_0(n)\rightarrow T^{\frac{n(n-1)}{2}};~A\mapsto \{\theta_1(u,A)_j^{(i)}\}_{1\le j\le i\le n-1}\]
with 
\[\theta_1(u,A)_j^{(i)}={\rm Arg}\left( (u_k-u_{k-1})^{\frac{\lambda^{(k)}_i(A)-A_{kk}}{2\pi\I }} (u_{k+1}-u_{k})^{\frac{A_{k+1,k+1}-\lambda^{(k)}_i(A)}{2\pi\I }}\right).\]
Therefore, the diffeomorphism $\psi(u)=\mathcal{X}_\theta^{-1}\circ {\rm Ad}_{g(u;\cdot )}=\mathcal{X}_{\theta^{-1}\cdot \theta_1(u)}$ is a also a transformation along the Thimm torus generated by the (pointwise product) map $\theta^{-1}\cdot \theta_1(u):\Herm_0(n)\rightarrow T^{\frac{n(n-1)}{2}}$. Then by the explicit expression of the components of the map $\theta:\Herm_0(n)\rightarrow T^{\frac{n(n-1)}{2}}$ given in \eqref{linearsys}, we obtain that 
\begin{cor}
For any $u\in\h_{\rm reg}(\mathbb{R})$, the diffeomorphism $\psi(u)$ coincides with the transformation $\mathcal{X}_{\phi}$, where $\phi=\theta^{-1}\cdot \theta_1(u)$ is given explicitly by (when restricts to the open dense subset $\Herm_0(n)$) \[\phi:\Herm_0(n)\rightarrow T^{\frac{n(n-1)}{2}};~A\mapsto \{\phi_j^{(i)}(A)\}_{1\le j\le i\le n-1}\] with the components
\begin{align}\nonumber
\phi_j^{(i)}(A)=&{\rm Arg}\left(\frac{ \prod_{v=1}^{i}\Gamma\left(1+\frac{\lambda^{(i)}_v-\lambda^{(i)}_j}{2\pi\I }\right)}{\prod_{v=1}^{i-1}\Gamma\left(1+\frac{\lambda^{(i-1)}_v-\lambda^{(i)}_j}{2\pi\I }\right)} \frac{\prod_{v=1}^{i}\Gamma\left(1+\frac{\lambda^{(i)}_v-\lambda^{(i)}_j}{2\pi\I }\right)}{\prod_{v=1}^{i+1}\Gamma\left(1+\frac{\lambda^{(i+1)}_v-\lambda^{(i)}_j}{2\pi\I }\right)}\right)+\pi\I \\
&+{\rm Arg}\left( (u_k-u_{k-1})^{\frac{\lambda^{(k)}_i(A)-A_{kk}}{2\pi\I }} (u_{k+1}-u_{k})^{\frac{A_{k+1,k+1}-\lambda^{(k)}_i(A)}{2\pi\I }}\right).
\end{align}
\end{cor}

\Addresses
\end{document}